\documentclass[11pt]{article}
\usepackage{amsfonts}

\newcommand{\oper}[2]{\newcommand{#1}{\mathop{\mathrm{#2}}\nolimits} }
\oper{\tr}{tr}
\oper{\adj}{adj}
\oper{\Div}{div}
\oper{\ad}{ad}
\oper{\Ad}{Ad}
\oper{\End}{End}
\oper{\Hom}{Hom}
\oper{\Aut}{Aut}
\oper{\SO}{SO}
\oper{\SP}{Sp}
\oper{\SU}{SU}
\oper{\GL}{GL}
\oper{\T}{T}
\oper{\U}{U}
\oper{\id}{I}
\oper{\ext}{Ext}
\oper{\rank}{rank}
\oper{\diag}{Diag}

\def\Box{\triangle^c}

\def\bR{{\bf R}}
\def\bH{{\bf H}}
           
\def\bC{{\bf C}}

\def\bn{{\bf n}}

\def\cp{{\bf C}{\bf P}}

\def\ci{{\cal I}}

\newtheorem{corollary}{Corollary}
\newtheorem{proposition}{Proposition}
\newtheorem{theorem}{Theorem}
\newtheorem{definition}{Definition}
\newtheorem{lemma}{Lemma}

\newcommand{\bproof}{\noindent{\it Proof: }}
\newcommand{\eproof}{\  q.~e.~d. \vspace{0.2in}}
\begin{document}

\title{Geometry of Hyper-K\"ahler Connections with Torsion}
\author{Gueo Grantcharov \thanks{
Address: Department of Mathematics, University of California at Riverside,
Riverside, CA 92521, U.S.A.. E-mail: geogran@math.ucr.edu. Partially supported
by Contract MM 809/1998 with Ministry of Education of Bulgaria and by Contract
238/1998 with the University of Sofia.} \and Yat Sun
Poon \thanks{
Address: Department of Mathematics, University of California at Riverside,
Riverside, CA 92521, U.S.A.. E-mail: ypoon@math.ucr.edu. Partially supported by
the IHES.}}
\maketitle


\noindent{{\bf   Abstract}\ \ The internal space of a N=4 supersymmetric model
with Wess-Zumino term has a connection with totally skew-symmetric
torsion and holonomy in $\SP(n)$. 
We study the mathematical background of this 
type of connections. In particular, we relate it to classical Hermitian geometry,
construct homogeneous as well as inhomogeneous examples, 
characterize it in terms of holomorphic data, 
develop its potential theory and reduction theory.}

\ 



\section{Introduction}

It has been known that the internal space for N=2
supersymmetric one-dimensional sigma model is a K\"ahler 
manifold \cite{Zumino}, and the
internal space for N=4 supersymmetric one-dimensional sigma model is a
hyper-K\"ahler manifold \cite{Curtright} \cite{Hitchin}. 
It means that there exists a torsion-free 
connection with holonomy in $\U(n)$ or $\SP(n)$ respectively
on the internal space.

It is also known for a fairy long time that when the
Wess-Zumino term is present in the sigma model, then the internal space has
linear connections with holonomy in $\U(n)$ or $\SP(n)$ depending on the
numbers of supersymmetry. However, the connection has torsion and the
torsion tensor is totally skew-symmetric \cite{Zumino} 
\cite{HP1} \cite{GPS}. 
The geometry of a connection with totally skew-symmetric
torsion and holonomy in $\U(n)$ is referred to KT-geometry by physicists.
When the holonomy is in $\SP(n)$, the geometry is referred to 
HKT-geometry.

If one ignores the metric and the connection of a HKT-geometry, 
the remaining object on the
manifold is a hypercomplex structure. The subject
of hypercomplex manifolds has been studied by many people since the publication
of \cite{Salamon2} and \cite{Boyer}. Considerable amount of information is 
known. It has a twistor correspondence \cite{Salamon2} \cite{PP1}. There are
homogeneous examples \cite{Joyce2}. There are inhomogeneous examples 
\cite{BGM} \cite{PP2}. 
There is a reduction construction modeled on symplectic reduction and hyper-K\"ahler
reduction \cite{Joyce1}. However, all these work focus on the hypercomplex
structure and the associated Obata connection which is a torsion-free connection
preserving the hypercomplex structure. What is not discussed in 
these work is hyper-Hermitian geometry.

On the other hand, Hermitian connections on almost Hermitian manifolds are studied 
rather thoroughly by Gauduchon \cite{Gauduchon}. He considered a subset of Hermitian
connections determined by the form of their torsion tensor, called canonical
connections.

Guided by physicists' work and based on the results on hypercomplex
manifolds, we review and further develop the theory of HKT-geometry. 
When  some of our observations are
re-interpretation of physicists' results, especially those in \cite{HP1} \cite{HP2}
and \cite{OP} \cite{Spindel}, some of the results
 in this paper are new.

In Section \ref{HKT Geometry},
we review the basic definitions of HKT-geometry along the line of classical 
Hermitian geometry  developed by Gauduchon \cite{Gauduchon}.  
Based on Joyce's construction of homogeneous hypercomplex manifolds \cite{Joyce2},
we review the construction of homogeneous HKT-geometry 
with respect to compact semi-simple Lie groups \cite{OP}. 

In Section \ref{Associated}, we find that a hyper-Hermitian manifold admits
HKT-connection if and only if for each complex structure, there is a holomorphic
(0,2)-form. This characterization easily implies that some hyper-Hermitian
structures are not HKT-structure. Furthermore,
when this characterization is given a twistorial interpretation, 
the associated object on the twistor space of the hypercomplex structure is
holomorphic with respect to a non-standard almost complex structure
 ${\cal J}_2$. This almost complex structure ${\cal J}_2$
 is first discussed by Eells and Salamon in a different context \cite{ES}. 
Since this almost complex structure 
is never integrable, we focus on the holomorphic (0,2)-forms. From this
perspective, we verify that there are HKT-structures on nilmanifolds, 
and that the twist of a HKT-manifold is again a HKT-manifold.

Based on results in
Section \ref{Associated}, we study potential theory for HKT-geometry
in Section \ref{Potential}. We shall see that local HKT-geometry is very flabby in
the sense that the existence of one generates many through a perturbation
of potential functions. In particular,
we show that hyper-K\"ahler potentials generate many HKT-potentials. 
The results in this section and Section \ref{Associated} allow us to construct
a large family of inhomogeneous HKT-structures on compact manifolds including
$S^1\times S^{4n+3}$.

Finally, a reduction theory based on hyper-K\"ahler reduction for HKT-geometry is
developed in Section \ref{Reduction}.

\section{Hyper-K\"ahler Geometry with Torsion}\label{HKT Geometry}

\subsection{K\"ahler Geometry with Torsion}

Let $M$ be a smooth manifold with Riemannian metric $g$ and an integrable
complex structure $J$. It is a Hermitian manifold if $g(JX, JY)=g(X, Y)$.
The  K\"ahler form $F$ is a type (1,1)-form defined by $F(X, Y)=g(JX, Y)$.

A linear connection $\nabla$ on $M$ is Hermitian if it preserves the metric $
g$ and the complex structure $J$. i.e., 
\[
\nabla g=0 \mbox{ and } \nabla J=0.
\]
Since the connection preserves the metric, it is uniquely determined by its
torsion tensor $T$. We shall also consider the following (3,0)-tensor 
\begin{equation}
c(X, Y, Z)=g(X, T(Y,Z)).
\end{equation}

Gauduchon  found that on any Hermitian manifold, the
collection of canonical Hermitian connections is an affine subspace of the space of
linear connections \cite{Gauduchon}. This affine subspace is at most one
dimensional. It is one
point if and only if the Hermitian manifold is K\"ahler, i.e., when the
K\"ahler form is closed, then the family of canonical
Hermitian connections collapses
to the Levi-Civita connection of the given metric. It is one-dimensional if
and only if the Hermitian manifold is non-K\"ahler. In the latter case,
there
are several distinguished Hermitian connections. For examples, Chern
connection and Lichnerwicz's \it first canonical connection \rm are
in this family. We are interested in another connection in this family.

Physicists find that the presence of the Wess-Zumino term in 
N=2 supersymmetry yields a Hermitian connection whose torsion $c$ is
totally skew-symmetric. In other words, $c$ is a 3-form. Such a connection
turns out to be another distinguished Hermitian connection \cite{Bismut} 
\cite{Gauduchon}. The geometry of such a connection is called by physicists
a KT-connection. Among some mathematicians, this connection is called the
Bismut connection. According to Gauduchon \cite{Gauduchon}, on any Hermitian
manifold, there exists a unique Hermitian connection whose torsion tensor $c$
is a 3-from. Moreover, the torsion form can be expressed in terms of the
complex structure and the K\"ahler form. Recall the following
definitions and convention \cite[Equations 2.8 and 2.15-2.17]{Besse}. For
any n-form $\omega$, when
\begin{equation}
(J\omega)(X_1, \dots, X_n):=(-1)^n\omega(JX_1, \dots, JX_n) \quad \mbox{ then }
\quad d^c\omega=(-1)^nJdJ\omega.
\end{equation}
And
\begin{equation}
\partial=\frac12(d+id^c)=\frac12(d+(-1)^niJdJ), \quad {\overline\partial}
=\frac12(d-id^c)=\frac12(d-(-1)^niJdJ).
\end{equation}
By \cite{Gauduchon}, the torsion 3-form of the Bismut connection is 
\begin{equation}
c(X,Y,Z)=-\frac12d^cF(X,Y,Z).
\end{equation}

\subsection{Hyper-K\"ahler Connection and HKT-Geometry}

Three complex structures $I_1, I_2$ and $I_3$ on $M$ form a
hypercomplex structure if 
\begin{equation}\label{quaternion}
I_1^2=I_2^2=I_3^2=-1, 
\quad 
\mbox{ and } 
\quad
I_1I_2=I_3=-I_2I_1.
\end{equation}
A triple of such complex structures is equivalent to the existence of a
2-sphere worth of integrable complex structures: 
\begin{equation}
\mathcal{I}=\{a_1I_1+a_2I_2+a_3I_3: a_1^2+a_2^2+a_3^2=1\}.
\end{equation}
When $g$ is a Riemannian metric on the manifold $M$ such that it is Hermitian
with respect to every complex structure in the hypercomplex structure, 
$(M, \mathcal{I}, g)$ is called a hyper-Hermitian manifold. Note that $g$ is
hyper-Hermitian if and
only if 
\begin{equation}
g(X, Y)=g(I_1X, I_1Y)=g(I_2X, I_2Y)=g(I_3X, I_3Y).
\end{equation}

On a hyper-Hermitian manifold, there are two natural torsion-free connections,
namely the Levi-Civita connection and the Obata connection. However,
in general 
the Levi-Civita connection does not preserve the hypercomplex structure
and the Obata connection does not preserve the metric.
We are interested in the following type of connections.

\begin{definition}
A linear connection $\nabla $ on a hyper-Hermitian manifold $
(M,\mathcal{I},g)$ is hyper-Hermitian if 
\begin{equation}
\nabla g=0,\quad \mbox{ and }\quad \nabla I_{1}=\nabla I_{2}=\nabla I_{3}=0.
\end{equation}
\end{definition}

\begin{definition}
A linear connection $\nabla $ on a hyper-Hermitian manifold $
(M,\mathcal{I},g)$ is hyper-K\"{a}hler if it is a hyper-Hermitian and its torsion
tensor is totally skew-symmetric. 
\end{definition}

A hyper-K\"ahler connection is referred to 
a HKT-connection in physics literature. The geometry of
this connection or this connection is also referred to a HKT-geometry.

Note that a HKT-connection is also the Bismut
connection for each complex structure in the given hypercomplex structure.
For the complex structures $\{I_1, I_2, I_3\}$, we consider
their corresponding K\"ahler forms $\{F_1, F_2, F_3\}$ and the complex
operators $\{d_1, d_2, d_3\}$ where $d_i = {d_i}^c$. Due to Gauduchon's
characterization of Bismut
connection, we have 

\begin{proposition}
A hyper-Hermitian manifold $(M,\mathcal{I},g)$ admits a
hyper-K\"{a}hler connection if and only if 
$d_{1}F_{1}=d_{2}F_{2}=d_{3}F_{3}$.
If it exists, it is unique. 
\end{proposition}

In view of the uniqueness, we say that $(M, \mathcal{I}, g)$
is a HKT-structure if it admits a hyper-K\"ahler connection. If the
hyper-K\"ahler connection is also torsion-free, then the HKT-structure is a
hyper-K\"ahler structure. 

\subsection{Homogeneous Examples}

Due to Joyce \cite{Joyce2}, there is a family of homogeneous
hypercomplex structures associated to any compact semi-simple Lie group. In
this section, we briefly review his construction and demonstrate, as
Opfermann and Papadopoulos did \cite{OP}, 
the existence of homogeneous HKT-connections. 

Let $G$ be a compact semi-simple Lie group. Let $U$ be a
maximal torus. Let $\mathfrak{g}$ and $\mathfrak{u}$ be their algebras.
Choose a system of ordered roots with respect to $\mathfrak{u}_\bC$. Let $
\alpha_1$ be a maximal positive root, and $\mathfrak{h}_1$ the dual space of 
$\alpha_1$. Let $\partial_1$ be the $\mathfrak{sp}(1)$-subalgebra of $
\mathfrak{g}$ such that its complexification is isomorphic to $\mathfrak{h}
_1\oplus\mathfrak{g}_{\alpha_1}\oplus\mathfrak{g}_{-\alpha_1}$ where $
\mathfrak{g}_{\alpha_1}$ and $\mathfrak{g}_{-\alpha_1}$ are the root spaces
for $\alpha_1$ and $-\alpha_1$ respectively. Let $\mathfrak{b}_1$ be the
centralizer of $\partial_1$. Then there is a vector subspace $\mathfrak{f}_1$
composed of root spaces such that $\mathfrak{g}=\mathfrak{b}
_1\oplus\partial_1\oplus \mathfrak{f}_1$. If $\mathfrak{b}_1$ is not
Abelian, Joyce applies this decomposition to it. By inductively searching
for $\mathfrak{sp}(1)$ subalgebras, he finds the following \cite[Lemma 4.1]
{Joyce2}. 

\begin{lemma}
The Lie algebra $\mathfrak{g}$ of a compact Lie group $G$
decomposes as 
\begin{equation}
\mathfrak{g}=\mathfrak{b}\oplus _{j=1}^{n}\mathfrak{\partial}_{j}\oplus
_{j=1}^{n}\mathfrak{f}_{j},  \label{decomposition}
\end{equation}
with the following properties. {\rm (1)} $\mathfrak{b}$ is Abelian and $
\mathfrak{\partial}_{j}$ is isomorphic to $\mathfrak{sp}(1)$. {\rm (2)}
$\mathfrak{b}\oplus _{j=1}^{n}\mathfrak{\partial}_{j}$ contains $\mathfrak{u}$. 
{\rm (3)} Set $\mathfrak{b}_{0}=\mathfrak{g}$, $\mathfrak{b}_{n}=\mathfrak{b}$
and $\mathfrak{b}_{k}=\mathfrak{b}\oplus _{j=k+1}^{n}\partial _{j}\oplus
_{j=k+1}^{n}\mathfrak{f}_{j}$. Then $[\mathfrak{b}_{k},\partial _{j}]=0$ for 
$k\geq j$. {\rm (4)} $[\mathfrak{\partial}_{l},\mathfrak{f}_{l}]\subset 
\mathfrak{f}_{l}$. {\rm (5)} The adjoint representation of $\mathfrak{\partial}_{l}
$ on $\mathfrak{f}_{l}$ is reducible to a direct sum of the irreducible
2-dimensional representations of $\mathfrak{sp}(1)$. 
\label{joyce decomposition} 
\end{lemma}

When the group $G$ is semi-simple, the Killing-Cartan form 
is a negative definite inner product on the vector space $\mathfrak{g}$. 

\begin{lemma}
The Joyce Decomposition of a compact semi-simple Lie algebra
is an orthogonal decomposition with respect to the Killing-Cartan form.
\end{lemma}
\bproof  Since Joyce Decomposition given as in 
(\ref{decomposition}) is inductively
defined, it suffices to prove that the decomposition 
\begin{equation}
\mathfrak{g}_\bC=\mathfrak{b}_1\oplus\partial_1\oplus\mathfrak{f}_1
\end{equation}
is orthogonal. Recall that 
\begin{eqnarray*}
\partial_1 &=& \langle \mathfrak{h}_1, X_{\alpha_1}, X_{-\alpha_1}\rangle, 
\mbox{ } \mathfrak{f}_1=\oplus_{\alpha_1\neq \alpha>0, \langle \alpha,
\alpha_1\rangle\neq 0} \mathfrak{g}_{\alpha}\oplus\mathfrak{g}_{-\alpha}, \\
\mathfrak{b}_1&=&\{h\in\mathfrak{u}_\bC: \alpha_1(h)=0\}
\oplus_{\alpha_1\neq\alpha>0, \langle \alpha, \alpha_1\rangle=0} \mathfrak{g}
_{\alpha}\oplus\mathfrak{g}_{-\alpha}.
\end{eqnarray*}
Since the Cartan subalgebra $\mathfrak{u}_\bC$ is orthogonal to any root
space, and it is an elementary fact that two root spaces $\mathfrak{g}
_{\alpha}$, $\mathfrak{g}_{\beta}$ are orthogonal whenever $\alpha\neq \pm
\beta$, $\mathfrak{f}_1$ is orthogonal to both $\mathfrak{b}_1$
and $\partial_1$. For the same reasons, $\partial_1$ is orthogonal to the
summand $\oplus_{\alpha>0, \langle \alpha, \alpha_1\rangle=0} \mathfrak{g}
_{\alpha}\oplus\mathfrak{g}_{-\alpha}$ in $\mathfrak{b}_1$, and $\mathfrak{b}
_1$ is orthogonal to the summand $\langle X_{\alpha_1}, X_{-\alpha_1}\rangle$
in $\partial_1$. Then $\partial_1$ is orthogonal to $\mathfrak{b}_1$
because for any element
$h$ in the Cartan subalgebra in $\mathfrak{b}_1$, $\langle
h, h_1\rangle=\alpha_1(h)=0$. \ q.~e.~d. \vspace{0.2in} 

Let $G$ be a compact semi-simple Lie group with rank $r$. Then 
\begin{equation}
(2n-r)\mathfrak{u}(1)\oplus\mathfrak{g}\cong {\mathbf{R}}^n\oplus_{j=1}^n
\mathfrak{\partial}_j\oplus_{j=1}^n\mathfrak{f}_j.  \label{vvv}
\end{equation}
At the tangent space of the identity element of $T^{2n-r}\times G$, i.e.\
the Lie algebra $(2n-r)\mathfrak{u}(1)\oplus \mathfrak{g}$, a hypercomplex
structure $\{I_1, I_2, I_3\}$ is defined as follows. Let $\{ E_1, \dots,
E_n\}$ be a basis for ${\mathbf{R}}^n$. Choose isomorphisms $\phi_j$ from $
\mathfrak{sp}(1)$,  the real vector space of imaginary quaternions, to $
\mathfrak{\partial}_j$. It gives a real linear identification from  the
quaternions $\mathbf{H}$ to $\langle E_j\rangle\oplus\mathfrak{\partial}_j$.
If $H_j$, $X_j$ and $Y_j$ forms a basis for $\mathfrak{\partial}_j$ such
that $[H_j, X_j]=2Y_j$ and $[H_j, Y_j]=-2X_j$, then 
\begin{equation}
I_1E_j=H_j, I_2E_j=X_j, I_3E_j=Y_j.
\end{equation}
Define the action of $I_a$ on $\mathfrak{f}_j$ by $I_a(v)=[v,
\phi_j(\iota_a)]$ where $\iota_1=i, \iota_2=j, \iota_3=k$. The complex
structures $\{I_1, I_2, I_3\}$ at the other points of the group $
T^{2n-r}\times G$ are obtained by left translations. These complex
structures are integrable and form a hypercomplex structure
\cite{Joyce2}.

\begin{lemma}
When $G$ is a compact semi-simple Lie group with rank $r$, there exists a
negative definite bilinear form $\hat{B}$ on the decomposition 
$(2n-r)\mathfrak{u}(1)\oplus \mathfrak{g}\cong {\mathbf{R}}^{n}\oplus
_{j=1}^{n}\mathfrak{\partial}_{j}\oplus _{j=1}^{n}\mathfrak{f}_{j}$
such that {\rm (1)}
 its restriction to $\mathfrak{g}$ is the Killing-Cartan form,
{\rm (2)} it is hyper-Hermitian respect to the hypercomplex structure, 
and {\rm (3)} the above decomposition is orthogonal. 
\end{lemma}
\bproof  In $\partial_j$, we choose an orthogonal
basis $\{H_j, X_j, Y_j\}$ such that $H_j$ is in the Cartan subalgebra and 
\begin{equation}
B(H_j, H_j)=B(X_j, X_j)=B(Y_j, Y_j)=-\lambda_j^2.
\end{equation}
On ${\mathbf{R}}^n=(2n-r)\mathfrak{u}(1)\oplus\mathfrak{b}$, choose $E_1,
\dots, E_n$ and extend the Killing-Cartan form so that 
\begin{equation}
{\hat B}(E_i, E_j)=-\delta_{ij}\lambda_j^2.
\end{equation}
It is now apparent that the extended Killing-Cartan form is hyper-Hermitian
with respect to $I_1, I_2$ and $I_3$. 

To show that the Killing-Cartan form is hyper-Hermitian on $
\oplus_{j=1}^n\mathfrak{f}_j$, it suffices to verify that the Killing-Cartan
form is hyper-Hermitian on $\mathfrak{f}_1$. It follows from the
fact that $B(X, [Y, Z])$ is totally skew-symmetric with respect in
$X, Y, Z$ and the Jacobi identity. \eproof

Let $g$ be the left-translation of the extended
Killing-Cartan form $-\hat B$. It is a bi-invariant metric on the manifold $
T^{2n-r}\times G$. The Levi-Civita connection $D$ is the
bi-invariant connection. Let $\nabla$ be the left-invariant connection. When 
$X$ and $Y$ are left-invariant vector fields 
\[
D_XY=\frac{1}{2}[X,Y], \mbox{ and } \nabla_XY=0. 
\]
Since the hypercomplex structure and the hyper-Hermitian metric are
left-invariant, the left-invariant connection is hyper-Hermitian.
The torsion tensor for the left-invariant connection is
$T(X,Y)=-[X, Y]$. The (3,0)-torsion tensor is 
\[
c(X,Y,Z)=-{\hat B}([X,Y],Z). 
\]
It is well-known that $c$ is a totally skew-symmetric 3-form. Therefore, 
the left-invariant connection is a  HKT-structure
on the group manifold $T^{2n-r}\times G$. 

It is apparent that if one extends the Killing-Cartan form
in an arbitrary way, then the resulting bi-invariant metric and
left-invariant hypercomplex structure cannot make a hyper-Hermitian
structure. 

The above construction can be generalized to homogeneous
spaces \cite{OP}.

\section{Characterization of HKT-Structures}\label{Associated}

In this section, we characterize HKT-structures in terms of the
existence of a holomorphic object with respect to any complex
structure in the hypercomplex structure. Through this characterization,
we shall find other examples of HKT-manifolds. Toward the end of
this section, we shall also
reinterpret the twistor theory for HKT-geometry
developed by Howe and Papadopoulos \cite{HP2}.
The results seem to indicate that the holomorphic characterization developed in the
next paragraph will serve all the purpose that one wants the twistor theory
of HKT-geometry to serve.

\subsection{Holomorphic Characterization}

\begin{proposition}\label{character}
Let $(M,\mathcal{I},g)$ be a hyper-Hermitian manifold and $
F_{a}$ be the K\"{a}hler form for $(I_{a},g)$. Then $(M,\mathcal{I},g)$ is a
HKT-structure if and only if $\partial _{1}(F_{2}+iF_{3})=0$;
 or equivalently 
${\overline{\partial }}_{1}(F_{2}-iF_{3})=0$.
\end{proposition}
\bproof  Since 
$\partial_1 (F_2 + iF_3) 
= \frac{1}{2}(dF_2 - d_1F_3) + \frac{i}{2}(d_1F_2 + dF_3)$,
it is identically zero if and only if $d_1F_2=-dF_3$,  and 
$dF_2=d_1F_3$. 

Note that
$F_2(I_1X, I_1Y) = g(I_2 I_1 X, I_1Y) = -g(I_2X, Y) =
-F_2(X, Y)$. It follows that
$d_1F_2 = (-1)^2I_1 d I_1 (F_2) = -I_1 dF_2$.
As $dF_2$ is a 3-form, for any $X, Y, Z$ tangent vectors, 
\begin{eqnarray*}
-I_1dF_2(X,Y,Z) &=& dF_2(I_1X, I_1Y, I_1Z)=dF_2(I_2I_3X, I_2I_3Y, I_2I_3Z) \\
&=& -I_2dF_2(I_3X, I_3Y, I_3Z)=I_3I_2dF_2(X, Y, Z).
\end{eqnarray*}
Since $F_2$ is type (1,1) with respect to $I_2$, $I_2F_2=F_2$. Then 
$d_1F_2 =-I_1dF_2=I_3I_2dF_2=I_3I_2dI_2F_2=I_3d_2F_2$. 
On the other hand, 
$-dF_3=I_3I_3dF_3=I_3I_3dI_3F_3=I_3d_3F_3$. 
Therefore, $d_2F_2=d_3F_3$ if and only if $d_1F_2=-dF_3$. Similarly, one can
prove that $d_2F_2=d_3F_3$ if and only if $d_1F_3 = dF_2$. It follows that 
$\partial_1(F_2+iF_3)=0$ if and only if $d_2F_2=d_3F_3$.
It is equivalent to 
$\nabla^2=\nabla^3$ 
where $\nabla^a$ is the Bismut connection of the Hermitian
structure $(M, I_a, g)$. Since $I_1=I_2I_3$, and $\nabla^2=\nabla^3$, $I_1$
is parallel with respect to $\nabla^2=\nabla^3$. By the uniqueness of Bismut
connection, $\nabla^1=\nabla^2=\nabla^3$. \eproof

On any hypercomplex manifold $(M, {\cal I})$, 
if $F_2-iF_3$ is a 2-form such that 
 $-F_{2}(I_{2}X,Y)=g(X,Y)$ is  positive definite and 
it is a non-holomorphic (0,2)-form with respect to $I_1$,
 then $(M, g, {\cal I})$ is a hyper-Hermitian
manifold but it is not a HKT-structure. For example, a conformal change
of a HKT-structure by a generic function gives a hyper-Hermitian structure
which is not a HKT-structure so long as the dimension of the underlying 
manifold is at least eight. On the other hand Proposition \ref{character}
implies that every four-dimensional hyper-Hermitian manifold is a 
HKT-structure, a fact also proven in \cite[Section 2.2]{GT}.

In the proof of Proposition \ref{character}, we also derive the
following \cite{HP2}. 

\begin{corollary}
Suppose $F_{1},F_{2}$ and $F_{3}$ are the K\"{a}hler forms of
a hyper-Hermitian structure. Then the hyper-Hermitian structure is a
HKT-structure if and only if 
\begin{equation}\label{difj}
d_{i}F_{j}=-{2}\delta _{ij}c-\epsilon _{ijk}dF_{k}.
\end{equation}
\end{corollary}

\begin{theorem}\label{holomorphic}
Let $(M,\mathcal{I})$ be a hypercomplex manifold and 
$F_{2}-iF_{3}$ be a {\rm (0,2)}-form with respect to $I_{1}$ such
that  ${\overline{\partial }}_{1}(F_{2}-iF_{3})=0$ or equivalently
${\partial }_{1}(F_{2}+iF_{3})=0$
and  $-F_{2}(I_{2}X,Y)=g(X,Y)$ is a positive definite
symmetric bilinear form. 
Then $(M,\mathcal{I},g)$ is a HKT-structure. 
\end{theorem}
\bproof In view of the last proposition, it suffices to
prove that the metric $g$ along with the given hypercomplex structure $
\mathcal{I}$ is hyper-Hermitian. 

Note that $F_2-iF_3$ is type (0,2) with respect to 
$I_1$. Since $X-iI_1X$ is a type (1,0)-vector with respect to $I_1$,
$(F_2-iF_3)(X-iI_1X, Y)=0$ for any vectors $X$ and $Y$.
It is equivalent to the identity $F_2(I_1X, Y) = -F_3(X, Y).$
Then 
\[
F_3(I_3X, Y) = -F_2(I_1I_3X, Y) = F_2(I_2X, Y)=-g(X, Y). 
\]
So $F_3(I_3X, I_3Y) = F_3(X, Y)$, and $g$ is Hermitian with respect to
$I_3$.
Since the metric $g$ is Hermitian with respect
to $I_2$ and $I_1=I_2I_3$, $g$ is also Hermitian with respect $I_1$. 
\eproof 

\subsection{HKT-Structures on Compact Nilmanifolds}

In this section, we apply the last theorem to construct a homogeneous
HKT-structure on some compact nilmanifolds.

Let $\{ X_1, ... , X_{2n},Y_1, ... , Y_{2n}, Z \}$  be a basis for 
${\mathbf R}^{4n + 1}$. Define commutators by:
$[X_i, Y_i] = Z$, and all others are zero. 
These commutators define on ${\mathbf R}^{4n + 1}$ the structure of
the \it Heisenberg Lie algebra \rm  ${\it h}_{2n}$.
Let ${\mathbf R}^3$ be the
3-dimensional Abelian algebra. 
The direct sum  $\bn = {\it h}_{2n} \oplus {\mathbf R}^3$ is a 2-step
nilpotent algebra whose center is four dimensional. 
Fix a basis $\{ E_1, E_2, E_3\}$ for 
${\mathbf R}^3$ and consider the following endomorphisms 
of ${\bf n}$ \cite{Dotti} :
\begin{eqnarray*}
I_1 &:& X_i \rightarrow Y_i,  Z \rightarrow E_1, E_2 \rightarrow
E_3;\\
I_2 &:& X_{2i+1} \rightarrow X_{2i},  Y_{2i-1} \rightarrow
Y_{2i},  Z \rightarrow E_2, E_1 \rightarrow E_3;
\\
I_1^2 &=& I_2^2=-{\rm identity}, \hspace{.5in} I_3 = I_1I_2.
\end{eqnarray*}
Clearly $I_1I_2 = -I_2I_1$. Moreover,
for $a = 1, 2, 3$ and $X, Y \in \bn $, 
$[I_aX, I_aY] = [X, Y]$ so $I_a$ are Abelian
complex structures on 
$\bn$ in the sense of \cite{Dotti} and in particular are
integrable.
It implies that $ \{ I_a : a =1, 2, 3 \} $ is a left invariant hypercomplex structure
on the simply connected Lie group $N$ whose algebra is $\bf n$. 
It is known that the complex structures $I_a$ on $\bf n$
satisfy:
\[
d( \Lambda^{1,0}_{I_a} \bn^* ) \in  \Lambda^{1,1}_{I_a} \bn^*
\]
where $\bn^* $ is the space of left invariant  1-forms on $N$ and 
$ \Lambda^{i,j}_{I_a} \bn^* $
is the $(i, j)$-component of ${\bf n}^* \otimes {\bf C}$
with respect to $I_a$ \cite{Dotti}. But then we have
$d( \Lambda^{2,0}_{I_a} \bn^* ) \in  \Lambda^{2,1}_{I_a} \bn^*$
and any left invariant (2,0)-form is $\partial_1$-closed. 
Now consider
the invariant metric on $N$ for which the basis $ \{ X_i , Y_i , Z, E_a \}$ is
orthonormal. Since it is compatible with the structures $I_a$ in view of
Theorem \ref{holomorphic} we obtain a left-invariant HKT-structure on $N$.
Noting that $N$ is isomorphic to the product $H_{2n} \times {\mathbf R}^3$
of the Heisenberg Lie group $H_{2n}$ and the Abelian group ${\mathbf R}^3$
we have:
 \begin{corollary}
Let ${\Gamma}$ be a cocompact lattice in the Heisenberg
group $H_{2n}$ and ${\bf Z}^3$ a lattice in ${\mathbf R}^3$.
The compact nilmanifold $(\Gamma \times {\bf Z}^3 ) \backslash N$
admits a HKT-structure.
\end{corollary}

\subsection{Twist of Hyper-K\"ahler Manifolds with Torsions}

Suppose that $(M, {\cal I})$ is a hypercomplex manifold, a $\U(1)$-instanton
$P$ is a principal $\U(1)$-bundle with a $\U(1)$-connection 1-form $\theta$
 such that its curvature 2-form is type-(1,1)
with respect to every complex structure in ${\cal I}$ \cite{CS} \cite{GP}.
Let $\Psi_M:\U(1)\to \Aut (M)$ be a group of hypercomplex automorphism,
and let $\Psi_P:\U(1)\to \Aut (P)$ be a lifting of $\Psi_M$. Let
$\Phi:\U(1)\to \Aut P$ be the principal $U(1)$-action on the bundle $P$,
and $\triangle (g)$ be the diagonal product $\Phi(g)\Psi_P(g)$ action on $P$.
A theorem of Joyce \cite[Theorem 2.2]{Joyce2} states that  
the quotient space $W=P/\triangle (\U(1))$
of the total space of $P$ with
respect to the diagonal action $\triangle$ is a hypercomplex manifold
whenever the vector fields generated by $\triangle (\U(1))$ are 
transversal
to the horizontal distribution of the connection $\theta$.
The quotient space $W$ is called a twist of the hypercomplex manifold
$M$.

Now suppose that $(M, \ci, g)$ is a HKT-structure, and $P$
is a $\U(1)$-instanton with connection form $\theta$. 
Suppose that $\Psi_M: \U(1)\to \Aut(M)$
is a group of hypercomplex isometry. Due to the uniqueness of HKT-structure,
$\Psi_M$ is a group of automorphism of the HKT-structure. 

\begin{corollary}
The twist manifold $W$ admits a HKT-structure.
\end{corollary}
\bproof
Let $\phi : P \rightarrow M$ and 
$\Delta  : P \rightarrow W$
be the projections from the instanton bundle P to $M$  and the twist $W$
respectively.  The connection $\theta$ defines a splitting of the 
tangent bundle of P into horizontal and vertical components: 
$TP = \cal{H} \oplus \cal{V}$ 
where $\cal{H} = \it{Ker} \theta$. We define 
endomorphisms $\tilde{I_a}$ on $TP$ as follows:
$\tilde{I_a} = 0$
on vertical directions, and
when $\tilde v$ is a horizontal lift of a tangent vector $v$ to $M$,
define $ \tilde{I_a} \tilde{v} = \widetilde {I_a v}$. 

Since the fibers of the projection $\Delta$ is
transversal to the horizontal distribution, for any tangent vector $\hat v$ to
$W$, there exists a horizontal vector $\tilde v$ such that 
$d\Delta {\tilde v}={\hat v}$.
Define ${\hat I}_a$ and $\hat g$  on $W$ by
${\hat I}_a {\hat v} = d\Delta(\tilde{I_a} \tilde{v})$ and
${\hat g}({\hat v},{\hat w}) = \tilde{g}( \tilde{v}, \tilde{w})$.
As the diagonal action is a group of hyper-holomorphic isometries, 
the almost complex structures ${\hat I}_a$ and metric $\hat g$ are
well-defined.

To verify that ${\hat I}_a$ are integrable complex structures on $W$, 
 we first observe that:
for horizontal vector fields $X$ and $Y$,
$d\Delta [X, Y] = [d\Delta X, d\Delta Y]$, 
 $d\phi [X, Y] = [d\phi X, d\phi Y]$ and
  $d\Delta {\tilde I}_a = {\hat I}_a d\Delta$, 
 $d\phi{\tilde I}_a = I_a d\phi$.
Through these relations, we establish the following
relations between
 Nijenhius tensors of $I_a$, ${\hat I}_a$ and $\tilde{I_a}$:
\[
d\Delta \tilde{N_a} (X, Y) = {\hat N}_a (d\Delta X, d\Delta Y)
\hspace{.2in}
\mbox{ and }
\hspace{.2in}
d\phi\tilde{N_a} (X, Y) = N_a (d\phi X, d\phi Y).
\]
The second identity implies that the horizontal part of 
$\tilde{N_a} (X, Y)$ vanishes because the complex structures
$I_a$ are integrable. With the first identity, it follows that
the Nijenhius tensor for ${\hat I}_a$ vanishes
if the  vertical part of $\tilde{N_a} (X, Y)$ also vanishes.
To calculate the vertical part, we have 
\begin{eqnarray*}
\theta (\tilde{N_a} (X, Y)) &=&
\frac14 \theta ([X, Y]+{\tilde I}_a[{\tilde I}_aX, Y]
+{\tilde I}_a[X, {\tilde I}_aY]-[{\tilde I}_aX, {\tilde I}_aY])\\
&=&
\frac14 \theta ([X, Y]-[{\tilde I}_aX, {\tilde I}_aY])
= \frac14(d\theta(X, Y) - d\theta(I_aX, I_aY)).
\end{eqnarray*}
Since $\theta$ is an instanton,  $d\theta(X, Y) - d\theta(I_aX, I_aY)=0$.
It follows that ${\hat I}_a$ are integrable.

To check that $\hat{g}$ is a HKT-metric, we first observe
that $d\Delta$ and $d\phi$ give rise to isomorphisms of $\Lambda^{(p,q)}M$,
$\Lambda^{(p,q)} \cal{H}$ and $\Lambda^{(p,q)} W$ when we fix the structures $I_1$, ${\hat
I}_1$ and ${\tilde I}_1$. Let the K\"ahler forms of 
the structures $I_a$ and ${\hat I}_a$ be denoted
by $F_a$ and ${\hat F}_a$
respectively. Now if $X, Y$ and $Z$ are sections of
${\cal H}^{(1,0)}$ then
\[
X(\Delta^* ({\hat F}_2 + i{\hat F}_3))(Y, Z) = X(\phi^* (F_2 +
iF_3))(Y, Z).
\]
Since $d\theta$ is type (1,1), $\theta([X, Y]) = d\theta(X, Y) = 0$. It
means that $[X, Y]$ is a section of ${\cal H}^{(1,0)}$. Therefore,
$\Delta^* ({\hat F}_2 + i{\hat F}_3)([X, Y], Z) 
= \phi^* (F_2 +iF_3)([X, Y], Z)$.
It follows that
\[
(\Delta^* d({\hat F}_2 + i{\hat F}_3))|_{\Lambda^{(3,0)} \cal{H}}
   =
(d\Delta^* ({\hat F}_2 + i{\hat F}_3))|_{\Lambda^{(3,0)} \cal{H}}
   = d\phi^* (F_2 +iF_3))|_{\Lambda^{(3,0)} \cal{H}} = 0.
\] 
Hence $d({\hat F}_2 + i{\hat F}_3)|_{\Lambda^{(3,0)}W}= 0$ 
and the corollary follows from Proposition
\ref{character}.
\eproof

\subsection{Twistor Theory of HKT-Geometry}

When $(M, \mathcal{I})$ is a 4n-dimensional hypercomplex
manifold, the smooth manifold $Z=M\times S^2$ admits an integrable complex
structure. It is defined as follows. For a unit vector ${{\vec{a}}}=(a_1,
a_2, a_3)\in\mathbf{R}^3$, let $I_{{\vec{a}}}$ be the complex structure $
a_1I_1+a_2I_2+a_3I_3$ in the hypercomplex structure $\mathcal{I}$. Let $J_{{
\vec{a}}}$ be the complex structure on $S^2$ defined by cross product in $
\mathbf{R}^3$: $J_{{\vec{a}}}{{\vec{w}}}={{\vec{a}}}\times{{\vec{w}}}$. Then
the complex structure on $Z=M\times S^2$ at the point $(x, {{\vec{a}}})$ is 
$\mathcal{J}_{(x, {{\vec{a}}})}=I_{{\vec{a}}}\oplus J_{{\vec{a}}}.$
It is well-known from twistor theory that this complex structure is
integrable \cite{Salamon2}. 
We shall have to consider a non-integrable almost complex
structure $\mathcal{J}_2=I\oplus (-J).$
Unless specified the otherwise, we discuss holomorphicity on $Z$ in terms of
the integrable complex structure $\mathcal{J}$. 

With respect to $\mathcal{J}$, the fibers of the projection $
\pi$ from $Z=M\times S^2$ onto its first factor are holomorphic curves with
genus zero. It can be proved that the holomorphic normal bundles are $
\oplus^{2n}\mathcal{O}(1)$. The antipodal map $\tau$ on the second factor is
an anti-holomorphic map on the twistor space $Z$  leaving the fibers of the
projection $\pi$ invariant. 

The projection $p$ onto the second smooth factor of $
Z=M\times S^2$ is a holomorphic map such that the inverse image of a point $
(a_1, a_2, a_3)$ is the manifold $M$ equipped with the complex structure $
a_1I_1+a_2I_2+a_3I_3$. If $\mathcal{D}$ is the sheaf of kernel of the
differential $dp$, then we have the exact sequence 
\begin{equation}
0\to \mathcal{D}\to \Theta_Z \stackrel{dp}{\longrightarrow} p^*\Theta_{
\mathbf{C}\mathbf{P}^1} \to 0.
\end{equation}
Real sections, i.e. $\tau$-invariant sections,
of the holomorphic projection $p$ are fibers of
the projection from $Z$ onto $M$. 

Twistor theory shows that there is a one-to-one
correspondence between hypercomplex manifold $(M, \mathcal{I})$ and its
twistor space $Z$ with the complex structure $\mathcal{J}$, the
anti-holomorphic map $\tau$, the holomorphic projection $p$ and the sections
of the projection $p$ with prescribed normal bundle \cite{PP1}. 

It is not surprising that when a hypercomplex manifold 
has a HKT-structure, there is an additional geometric structure on the twistor
space. The following theorem is essentially 
developed in \cite{HP2}. 

\begin{theorem}
Let $(M,\mathcal{I},g)$ be a 4n-dimensional HKT-structure.
Then the twistor space $Z$ is a complex manifold such that 
\begin{enumerate}
\item  the fibers of the projection $\pi :Z\rightarrow M$
are rational curves with holomorphic normal bundle 
$\oplus ^{2n}\mathcal{O}(1)$, 
\item  there is a holomorphic projection $p:Z\rightarrow 
\mathbf{C}\mathbf{P}^{1}$ such that the fibers are the manifold $M$ equipped
with complex structures of the hypercomplex structure $\mathcal{I}$, 
\item  there is a ${\cal J}_{2}$-holomorphic section of $\wedge
^{(0,2)}\mathcal{D}\otimes p^{\ast }{\overline{\Theta }}_{\mathbf{C}\mathbf{P
}^{1}}$ defining a positive definite (0,2)-form on each fiber, 
\item  there is an anti-holomorphic map $\tau $ compatible
with 1, 2 and 3 and inducing the antipodal map on $\mathbf{C}\mathbf{P}
^{1}$. 
\end{enumerate}

Conversely, if $Z$ is a complex manifold with a
non-integrable almost complex structure $J_{2}$ with the above four
properties, then the parameter space of real sections of the projection $p$
is a 4n-dimensional manifold $M$ with a natural HKT-structure for which $Z$
is the twistor space. 
\end{theorem}
\bproof Given a HKT-structure, then only part {\it 3} in
the first half of this theorem is a new observation. It is a generalization
of Theorem \ref{holomorphic}. Through the stereographic projection, 
\begin{equation}
\zeta \mapsto {{\vec{a}}}
=\frac{1}{1+|\zeta |^{2}}
(1-|\zeta |^{2}, -i(\zeta -{\overline{\zeta }}), -(\zeta +{\overline{\zeta }}))
\end{equation}
$\zeta $ is a complex coordinate of the Riemann sphere. Note that 
\[
\frac{1}{1+|\zeta |^{2}}\left( 
\begin{array}{ccc}
1-|\zeta |^{2} 
	& i(\zeta -{\overline{\zeta }}) 
		& \zeta +{\overline{\zeta }}
\\ 
-i(\zeta -{\overline{\zeta }}) 
	& 1+\frac{1}{2}(\zeta ^{2}+{\overline{\zeta }}^{2}) 
		& -\frac{i}{2}(\zeta ^{2}-{\overline{\zeta }}^{2}) \\ 
-(\zeta +{\overline{\zeta }}) 
	& -\frac{i}{2}(\zeta ^{2}-{\overline{\zeta }}^{2}) 
		& 1-\frac{1}{2}(\zeta ^{2}+{\overline{\zeta }}^{2})
\end{array}
\right) 
\]
is a special orthogonal matrix. Let $\vec{b}$ and $\vec{c}$ be the second
and third column vectors respectively. Consider the complex structure 
\[
I_{\vec{a}}=\frac{1}{1+|\zeta |^{2}}\left( 
(1-|\zeta |^{2})I_{1}-i(\zeta -{\overline{\zeta }})I_{2}
-(\zeta +{\overline{\zeta }}) I_{3}\right).
\]
According to Theorem \ref{holomorphic}, the 2-form
\begin{equation}
F_{\vec{b}}-iF_{\vec{c}} 
=\frac{1}{1+|\zeta |^{2}}
\left( 
(F_{2}-iF_{3})-2i{\overline{\zeta }}F_{1}
+{\overline{\zeta }}^{2}(F_{2}+iF_{3})
\right)
\end{equation}
is holomorphic with respect to 
$I_{\vec{a}}$. 

Due to the integrability of the complex structure
 $I_{\vec{a}}$,  $d_{\vec{a}}$ is linear in $\vec{a}$. 
Therefore,
\begin{equation}\label{da}
d_{\vec{a}}=\frac{1}{1+|\zeta |^{2}}
\left( 
(1-|\zeta |^{2})d_{1}
-i(\zeta -{\overline{\zeta }})d_{2}
-(\zeta +{\overline{\zeta }}) d_{3}
\right).
\end{equation}

Note that ${\overline{\zeta }}$ is holomorphic with respect to the
almost complex structure ${\cal J}_{2}$. More precisely, consider the $\overline{
\partial }$-operator with respect to the almost complex structure ${\cal J}_{2}$:
on n-forms, it is
\begin{equation}
{\overline{\delta }}=\frac{1}{2}(d-i(-1)^{n}{\cal J}_{2}d{\cal J}_{2}),
\end{equation}
then ${\cal J}_{2}d{\overline{\zeta }}=i{\overline{\zeta }}$,
and  ${\overline{\delta }}
{\overline{\zeta }}=0$. It follows that at $(x,{\vec{a}})$ on $Z=M\times
S^{2}$, 
\begin{eqnarray*}
&&{\overline{\delta }}
\left( 
-2i{\overline{\zeta }}F_{1}
+(1+{\overline{\zeta }}^{2})F_{2}
-i(1-{\overline{\zeta }}^{2})F_{3}
\right)  
\\
&=&
-2i{\overline{\zeta }}{\overline{\delta }}F_{1}
+(1+{\overline{\zeta }}^{2}){\overline{\delta }}F_{2}
-i(1-{\overline{\zeta }}^{2}){\overline{\delta }}F_{3} 
=
-2i{\overline{\zeta }}{\overline{\partial }}_{\vec{a}}F_{1}
+(1+{\overline{\zeta }}^{2}){\overline{\partial }}_{\vec{a}}F_{2}
-i(1-{\overline{\zeta }}^{2}){\overline{\partial }}_{\vec{a}}F_{3}
\\
&=&
\frac12\left(
-2i{\overline{\zeta }}{d}F_{1}
+(1+{\overline{\zeta }}^{2}){d}F_{2}
-i(1-{\overline{\zeta }}^{2}){d}F_{3}
\right)\\
&&
-\frac{i}{2}\left(
-2i{\overline{\zeta }}{d}_{\vec{a}}F_{1}
+(1+{\overline{\zeta }}^{2}){d}_{\vec{a}}F_{2}
-i(1-{\overline{\zeta }}^{2}){d}_{\vec{a}}F_{3}
   \right)
\end{eqnarray*}
Now (\ref{da}) and (\ref{difj}) together imply that
 the twisted 2-form
$(F_{2}-iF_{3})-2i{\overline{\zeta }}F_{1}
+{\overline{\zeta }}^{2}(F_{2}+iF_{3})$
is closed with respect to $\overline\delta$. Therefore,
it is a ${\cal J}_2$-holomorphic section. 

Since $\zeta$ is a holomorphic coordinate on $S^2$, the homogeneity
shows that this section is twisted by $\overline{{\cal O}(2)}$.

The inverse construction is a consequence of the inverse construction
of hypercomplex manifold \cite{PP1} and Theorem \ref{holomorphic}.
\eproof

As the almost complex structure ${\cal J}_2$ is never integrable
\cite{ES}, twistor theory loses substantial power of holomorphic
geometry when we study HKT-structure. Therefore, we focus on 
application of Theorem \ref{holomorphic}.

\section{Potential Theory}\label{Potential}

Theorem \ref{holomorphic} shows that the form $F_2+iF_3$ is
a $\partial_1$-closed
(2,0)-form on a HKT-manifold. It is natural to consider a  differential form $\beta_1$
as potential 1-form for $F_2+iF_3$ if $\partial_1\beta_1=F_2+iF_3$. A priori,
the 1-form $\beta_1$ depends on the choice of the complex structure $I_1$.
The potential 1-form for $F_3+iF_1$, if it exists, depends on $I_2$, and so on.
In this section, 
we seek a function that generates all K\"ahler forms.

\subsection{Potential Functions}

A function $\mu$ is a  potential function for a hyper-K\"ahler manifold 
$(M, \mathcal{I}, g)$ if the K\"ahler forms $F_a$ are equal to $dd_a\mu$.
Since $d_a=(-1)^nI_adI_a$ on n-forms, $d_a\mu=I_ad\mu.$ Therefore, 
\[
d_1d_2\mu=d_1I_2d\mu=-I_1dI_1I_2d\mu=-I_1dI_3d\mu=-I_1dd_3\mu=-I_1\Omega_3
=\Omega_3=dd_3\mu. 
\]
 Now we generalize this concept
to HKT-manifolds.

\begin{definition}
Let $(M,\mathcal{I},g)$ be a HKT-structure with K\"{a}hler
forms $F_{1},F_{2}$ and $F_{3}$. A possibly locally defined function $\mu$
is a potential function for the HKT-structure if 
\begin{equation}\label{f123}
F_{1}=\frac{1}{2}(dd_{1}+d_{2}d_{3})\mu ,\quad F_{2}=\frac{1}{2}
(dd_{2}+d_{3}d_{1})\mu ,\quad F_{3}=\frac{1}{2}(dd_{3}+d_{1}d_{2})\mu .
\end{equation}
\end{definition}

Due to the identities
$dd_a+d_ad=0$  and $d_ad_b+d_bd_a=0$, 
$\mu$ is a potential function if and only if 
\[
F_{\vec{a}}=\frac{1}{2}(dd_{\vec{a}}+d_{\vec{b}}d_{\vec{c}})\mu,
\]
when ${\vec{a}}={\vec{b}}\times{\vec{c}}$ and $F_{\vec{a}}$
is the K\"ahler form for the complex structure $I_{\vec{a}}=a_1I_1+a_2I_2+a_3I_3$.
Moreover, the torsion
3-form is given by
 $d_1F_1=d_2F_2=d_3F_3=\frac12d_1d_2d_3\mu$. Furthermore, since 
$\partial_a=\frac{1}{2}(d+id_a)$ and 
${\overline\partial}_a=\frac{1}{2}(d-id_a)$, 
\begin{equation}
F_2+iF_3=\frac12 (dd_2+idd_3+id_1d_2-d_1d_3)\mu=
2\partial_1I_2{\overline\partial}_1\mu.
\end{equation}

Conversely, if a function
 $\mu$  satisfies the above identity, it
satisfies
the last two identities in (\ref{f123}). Since the metric is hyper-Hermitian, 
for any vectors $X$ and $Y$, $F_1(X, Y)=F_2(I_3X, Y)$. Through the integrability of
the complex structures $I_1, I_2, I_3$, the quaternion identities (\ref{quaternion})
 and the last two identities in (\ref{f123}),
one derives the first identity in (\ref{f123}). 
Therefore, we have the following theorem which justifies our definition for potential
functions.
\begin{theorem} 
Let $(M,\mathcal{I},g)$ be a HKT-structure with K\"{a}hler
form $F_{1},F_{2}$ and $F_{3}$. A possibly locally defined function $\mu $
is a potential function for the HKT-structure if 
\begin{equation}
F_2+iF_3=2\partial_1I_2{\overline\partial}_1\mu.
\end{equation}
\end{theorem}

In this context, a HKT-structure is hyper-K\"ahler if and
only if the potential function satisfies the following identities. 
\begin{equation}
dd_1\mu=d_2d_3\mu, \quad dd_2\mu=d_3d_2\mu, \quad dd_3\mu=d_1d_2\mu.
\end{equation}

\noindent{\bf Remark:}
 As in the K\"ahler case, compact manifolds do not admit globally 
defined HKT potential.
To verify, let $f$ be a potential function and $g$ be the corresponding induced
metric. Define the complex Laplacian of $f$ with respect to $g$:
$$
\overline{\partial}^* \overline{\partial} f = {\Box} f = g(dd_1f, F_1)
$$
Then 
$0 \leq 2g(F_1, F_1) = g(dd_1 f + d_2d_3 f, F_1) = 2 {\Box} f,$
because 
$$
g(d_2d_3 f, F_ 1) = g(-I_2dd_1 f, F_1) = -g(dd_1 f, I_2F_1) = g(dd_1 f, F_1) =
{\Box}f.
$$
Now the remark follows from the standard arguments involving maximum
principle for second order elliptic differential equation just like in the K{\"a}hler
case since  $\Box f$ does not have zero-order terms.

\

\noindent{\bf Remark:}
If we introduce the following 
quaternionic operators acting on quaternionic valued forms on the left:
$\partial^H=d+id_1+jd_2+kd_3$, and
${\overline\partial}^H=d-id_1-jd_2-kd_3$,
then a real-valued function $\mu$ is a HKT-potential if 
$\partial^H{\overline\partial}^H\mu=-2iF_1-2jF_2-2kF_3$.

If we identify $\bH^n$ with $\bC^{2n}$, we deduce from Theorem \ref{holomorphic} that any 
pluri-subharmonic function in domain $\bC^{2n}$ is an HKT-potential. The converse however is
wrong. As we shall see in Example \ref{hopf} the function
$\log(|z|^2+
|w|^2)$ is a 
 HKT potential in $\bC^{2n}\backslash\{ 0\}$ but is not
pluri-subharmonic. 

\

\noindent{\bf Remark:} Given a HKT-metric $g$ with K\"ahler forms $F_1$, $F_2$ and $F_3$,
for any real-valued function $\mu$ we consider
\[
{\hat F}_2+i{\hat F}_3=F_2+iF_3+\partial_1I_2{\overline\partial}_1\mu.
\]
According to Theorem \ref{holomorphic} and other results in this section, whenever
the form ${\hat g}(X, Y):=-{\hat F}_2(I_2X,Y)$ is positive definite, we obtain a new
HKT-metric with respect to the old hypercomplex structure.

\subsection{HKT-Potentials Generated by Hyper-K\"ahler Potentials}

Let $(M, \mathcal{I}, g)$ be a hyper-K\"ahler manifold with
hyper-K\"ahler potential $\mu$. The K\"ahler forms are given by
$\Omega_a=dd_a\mu$. We
consider HKT-structures generated by potential functions through $\mu$. 

\begin{theorem}\label{modification}
Suppose $(M, \mathcal{I}, g)$ is a hyper-K\"ahler manifold
with hyper-K\"ahler potential $\mu$. For any smooth function $f$ of one
variable, let $U$ be the open subset of $M$ on which $\mu$ is defined and 
\begin{equation}\label{Inequality}
f^{\prime}(\mu )+\frac{1}{4}f^{\prime\prime}(\mu )|\nabla\mu |^2>0.
\end{equation}
Define a symmetric bilinear form $\hat g$ by 
\begin{equation}\label{ghat}
{\hat g}=f^{\prime}(\mu )g+\frac{1}{4}f^{\prime\prime}(\mu ) (d\mu\otimes
d\mu+I_1d\mu\otimes I_1d\mu+ I_2d\mu\otimes I_2d\mu + I_3d\mu\otimes
I_3d\mu).
\end{equation}
Then $(U, \mathcal{I}, {\hat g})$ is a HKT-structure with $f(\mu)$ as its
potential. 
\end{theorem}
\bproof Since $\mu$ is a hyper-K\"ahler potential
for the metric $g$, 
$\Omega_2+i\Omega_3=2\partial_1I_2{\overline\partial}_1\mu.$
It follows that 
\begin{eqnarray*}
2\partial_1I_2{\overline\partial}_1f &=& 2\partial_1f^{\prime}(\mu)I_2{
\overline\partial}_1\mu 
= 2f^{\prime}(\mu)\partial_1I_2{\overline\partial}_1\mu
+2f^{\prime\prime}(\mu )\partial_1\mu\wedge I_2{\overline\partial}_1\mu \\
&=& f^{\prime}(\mu )(\Omega_2+i\Omega_3)+\frac12 f^{\prime\prime}(\mu)
(d\mu+id_1\mu )\wedge (I_2d\mu-iI_2d_1\mu ).
\end{eqnarray*}
When $F_2$ and $F_3$ are the real and imaginary part of $2\partial_1I_2{
\overline\partial}_1f $ respectively, then 
\begin{equation}
F_2=f^{\prime}(\mu )\Omega_2+\frac12f^{\prime\prime}(\mu )(d\mu\wedge
I_2d\mu+d_1\mu\wedge I_2d_1\mu ).
\end{equation}
It is now straight forward to verify that
$-{\hat F}_2(I_2X, Y)= {\hat g}(X, Y)$.
Therefore, $\hat g$ together with given hypercomplex structure defines a
HKT-structure with the function $f$ as its potential so long as $\hat g$ is
positive definite. 

Since $g$ is hyper-Hermitian, the vector fields $
Y_0=\nabla\mu$ and $Y_a=I_a\nabla\mu$ are mutually orthogonal with equal
length. At any point where $Y_0$ is not the zero vector, we extend $\{Y_0,
Y_1, Y_2, Y_3\}$ to an orthonormal frame with respect to the hyper-K\"ahler
metric $g$. Any vector $X$ can be written as 
$X=a_0Y_0+a_1Y_1+a_2Y_2+a_3Y_3+X^\perp$
where $X^\perp$ is in the orthogonal complement of $\{Y_0, Y_1, Y_2, Y_3\}$.
Note that 
\[
d\mu (X^\perp) = g(\nabla\mu, X^\perp )=0, \mbox{ and } I_ad\mu
(X^\perp)=-g(\nabla\mu, I_aX^\perp )=g(I_a\nabla\mu, X^\perp )=0.
\]
Also, for $1\leq a\neq b\leq 3$, 
\begin{eqnarray*}
d\mu (Y_a) &=& g(\nabla\mu, I_a\nabla\mu)=0, \quad d\mu (Y_0)=|\nabla\mu|^2,
\\
\ I_bd\mu (Y_a) &=& -g(\nabla\mu, I_bI_a\nabla\mu )=0, \quad I_ad\mu
(Y_a)=-g(\nabla\mu, I_a^2\nabla\mu )=|\nabla\mu|^2.
\end{eqnarray*}
Then 
\[
{\hat g}(X, X) 
=f^{\prime}(\mu )(\sum_{\ell=0}^3a_\ell^2)|\nabla\mu|^2 +
\frac{f^{\prime\prime}(\mu )}{4}(\sum_{\ell=0}^3a_\ell^2)|\nabla\mu|^4 
=(f^{\prime}(\mu )+\frac{f^{\prime\prime}(\mu )}{4}|\nabla\mu
|^2)(\sum_{\ell=0}^3a_\ell^2)|\nabla\mu|^2.
\]
Therefore, $\hat g$ is positive definite on the open set defined by the
inequality (\ref{Inequality}).
\eproof

Note that for any positive integer $m$, $f(\mu)=\mu^m$ satisfies
 (\ref{Inequality}) whenever $\mu$ is positive. So does 
$f(\mu )=e^{\mu}$. Therefore, if $g$ is a hyper-K\"ahler metric
with a positive potential function $\mu$, the following metrics
are HKT-metrics.
\begin{eqnarray*}
g_m &=&m\mu^{m-2}(\mu g+\frac{m-1}{4} 
(d\mu\otimes
d\mu+I_1d\mu\otimes I_1d\mu+ I_2d\mu\otimes I_2d\mu + I_3d\mu\otimes
I_3d\mu)),\\
g_\infty &=& e^\mu(g+\frac{1}{4}
(d\mu\otimes
d\mu+I_1d\mu\otimes I_1d\mu+ I_2d\mu\otimes I_2d\mu + I_3d\mu\otimes
I_3d\mu)).
\end{eqnarray*}

\subsection{Inhomogeneous HKT-Structures on $S^1\times
S^{4n-3}$}\label{hopf}

On the complex vector space $(\mathbf{C}^n\oplus 
\mathbf{C}^n)\backslash\{0\}$, let $(z_\alpha, w_\alpha)$, $1\leq\alpha\leq n
$, be its coordinates.  We define a hypercomplex
structure to contain this complex structure as follow. 
\[
\begin{array}{cccc}
I_1dz_\alpha =-idz_\alpha, & I_1dw_\alpha=-idw_\alpha, & I_1d{\overline z}
_\alpha = id{\overline z}_\alpha, & I_1d{\overline w}_\alpha=id{\overline w}
_\alpha. \\ 
I_2dz_\alpha=d{\overline w}_\alpha, & I_2dw_\alpha=-d{\overline z}_\alpha, & 
I_2d{\overline z}_\alpha = d{w}_\alpha, & I_2d{\overline w}_\alpha=-d{z}
_\alpha. \\ 
I_3dz_\alpha=id{\overline w}_\alpha, & I_3dw_\alpha=-id{\overline z}_\alpha,
& I_3d{\overline z}_\alpha = -id{w}_\alpha, & I_3d{\overline w}_\alpha=id{z}
_\alpha.
\end{array}
\]
The function $\mu=\frac12(|z|^2+|w|^2)$ is the hyper-K\"ahler potential for
the standard Euclidean metric:
\begin{equation}
g = \frac{1}{2}(dz_{\alpha}\otimes d{\overline z}_\alpha +d{\overline z}
_\alpha \otimes dz_{\alpha} +dw_{\alpha}\otimes d{\overline w}_\alpha +d{
\overline w}_\alpha \otimes dw_{\alpha}).
\end{equation}

Since $|\nabla\mu|^2=2\mu$, the function $f(\mu)=\ln\mu$ satisfies the
inequality (\ref{Inequality}) on $\bC^{2n}\backslash\{0\}$.  
By Theorem \ref{modification}, $\ln\mu$ is the HKT-potential for a HKT-metric  
 $\hat g$ on $\bC^{2n}\backslash\{0\}$.

Next for any real number $r$, with $0<r<1$, and $\theta_1,
\dots, \theta_n$ modulo $2\pi$, we consider the integer group $\langle
r\rangle$ generated by the following action on $(\mathbf{C}^n\oplus \mathbf{C
}^n)\backslash\{0\}$. 
\begin{equation}
(z_\alpha, w_\alpha)\mapsto (re^{i\theta_\alpha}z_{\alpha},
re^{-i\theta_\alpha} w_\alpha).
\end{equation}
One can check that the group $\langle r\rangle$ is a group of hypercomplex
transformations. As observed in \cite{PP2}, the quotient space of $(\mathbf{C}
^n\oplus \mathbf{C}^n)\backslash\{0\}$ with respect to $\langle r\rangle$ is
the manifold $S^1\times S^{4n-1}=S^1\times \SP(n)/\SP(n-1)$. Since the group 
$\langle r\rangle$ is also a group of isometries with respect to the
HKT-metric $\hat g$ determined by $f(\mu)=\ln \mu$, 
the HKT-structure descends from 
$(\mathbf{C}^n\oplus \mathbf{C}^n)\backslash\{0\}$ 
to a HKT-structure on $S^1\times
S^{4n-1}$. Since the hypercomplex structures on $S^1\times S^{4n-1}$ are
parametrized by $(r, \theta_1, \dots, \theta_n)$ and a generic hypercomplex
structure in this family is inhomogeneous \cite{PP2}, we obtain a family of
inhomogeneous HKT-structures on the manifold $S^1\times S^{4n-1}$. 

\begin{theorem} Every hypercomplex deformation of the homogeneous
hypercomplex structure on $S^1\times S^{4n-1}$ admits a HKT-metric.
\end{theorem}

Furthermore ${\hat F}_2+i{\hat F}_3=2\partial_1I_2{
\overline\partial}_1\mu$ descends to $S^1\times S^{4n-1}$. However, the
function $\mu$ does not descend to $S^1\times S^{4n-1}$. Therefore, this
(2,0)-form has a potential form $I_2{\overline\partial}_1\mu$ but
not a globally defined potential function. 

\subsection{Associated Bundles of Quaternionic K\"ahler
Manifolds}

When $M$ is a quaternionic K\"ahler manifold, i.e. the
holonomy of the Riemannian metric is contained in the group $\SP(n)\cdot\SP
(1)$, the representation of $\SP(1)$ on quaternions $\mathbf{H}$ defines an
associated fiber bundle $\mathcal{U}(M)$ over the smooth manifold $M$ with $
\mathbf{H}\backslash\{0\}/\mathbf{Z}_2$ as fiber. Swann finds that there is
a hyper-K\"ahler metric $g$ on $\mathcal{U}(M)$ whose potential function $\mu
$ is the length of the radius coordinate vector field along each fiber \cite{Swann}. 
As in the last example, $\ln\mu$ is the potential function of a HKT-structure with
metric $\hat g$. 

Again, the metric $\hat g$ and the hypercomplex structure
are both invariant of fiberwise real scalar multiplication. Therefore, the
HKT-structure with metric $\hat g$ descends to the compact quotients defined by
integer groups generated by fiberwise real scalar multiplications. 

\section{Reduction}\label{Reduction}

First of all, we recall the construction of hypercomplex reduction
developed by Joyce \cite{Joyce1}. Let $G$ be a
compact group of hypercomplex automorphisms on $M$. Denote the algebra of
hyper-holomorphic vector fields by $\mathfrak{g}$. Suppose that $\nu =(\nu
_{1},\nu _{2},\nu _{3}):M\longrightarrow \bR^{3}\otimes \mathfrak{g}$ is a 
$G$-equivariant map satisfying the following two conditions. The Cauchy-Riemann
condition: 
$I_{1}d\nu _{1}=I_{2}d\nu _{2}=I_{3}d\nu _{3}$,
and the transversality condition: 
$I_a d\nu_a(X)\neq 0$ for all $X\in \mathfrak{g}$.
 Any map satisfying these conditions is 
called a $G$-moment map. Given a point $\zeta =(\zeta
_{1},\zeta _{2},\zeta _{3})$ in $\bR^{3}\otimes \mathfrak{g}$, 
denote the level set $\nu ^{-1}(\zeta )$ by $P$. 
Since the map $\nu $ is $G$-equivariant, level sets are invariant if the
group $G$ is Abelian or if the point $\varsigma $ is invariant. Assuming
that the level set $P$ is invariant, and the action of $G$ on $P$ is free,
then the quotient space $N=P/G$ is a smooth manifold. 

Joyce proved that
the quotient space $N=P/G$ inherits a natural hypercomplex
structure \cite{Joyce1}. 
His construction runs as follows.
For each point $m$ in the space $P$, its tangent space is 
\[
T_{m}P=\{t\in T_{m}M:d\nu _{1}(t)=d\nu _{2}(t)=d\nu _{3}(t)=0\}.
\]
Consider the vector subspace 
\[
U_{m}=\{t\in T_{m}P:I_{1}d\nu _{1}(t)=I_{2}d\nu _{2}(t)=I_{3}d\nu
_{3}(t)=0\}.
\]
Due to the transversality condition, this space is transversal to the
vectors generated by elements in $\mathfrak{g}.$ Due to the Cauchy-Riemann
condition, this space is a vector subspace of $T_{m}P$ with co-dimension $
\dim \mathfrak{g}$, and hence it is a vector subspace of $T_{m}M$ with
co-dimension $4\dim \mathfrak{g}$. The same condition implies that, as a
subbundle of $TM_{|P}$, $U$ is closed under $I_a$. We call the
distribution $U$ the hypercomplex distribution of the map $\nu $. Let 
$\pi :P\to N$
be the quotient map. For any tangent vector $v$ at $\pi (m)$, there exists a
unique element $\widetilde{v}$ in $U_{m}$ such that $d\pi (\widetilde{v})=v$. 
The hypercomplex structure on $N$ is defined by 
\begin{equation}
I_av=d\pi (I_a\widetilde{v}),\quad \mbox{ i.e. }\quad \widetilde{I_av}
=I_a\widetilde{v}.
\end{equation}

\begin{theorem}\label{reduction}
Let $(M, \mathcal{I}, g)$ be a HKT-manifold. Suppose that $G$ is compact group
of hypercomplex isometries. Suppose that $\nu $ is a $G$-moment map such
that along the invariant level set $P=\nu ^{-1}(\zeta )$, the hypercomplex
distribution $U$ is orthogonal to the Killing vector fields generated by the
group $G$, then the quotient space $N=P/G$ inherits a natural
HKT-structure. 
\end{theorem}
\bproof
Under the condition of this theorem, the hypercomplex
distribution along the level set $P$ is identical to the orthogonal
distribution 
\[
H_{m}=\{t\in T_{m}P:g(t,X)=0, X\in {\mathfrak{g}}\}.
\]
Now, we define a metric structure $h$ at $T_{\pi (m)}N$ as follows. For $
v,w\in T_{\pi (m)}N,$ 
\begin{equation}
h_{\pi (m)}(v,w)=g_{m}(\widetilde{v},\widetilde{w}).
\end{equation}
It is obvious that this metric on $N$ is hyper-Hermitian. To find the
hyper-K\"{a}hler connection $D$ on the quotient space $N$, let $v$ and $w$ be
locally defined vector fields on the manifold $N$. They lift uniquely to 
$G$-invariant sections $\widetilde{v}$ and $\widetilde{w}$ of the bundle $U$.
As $U$ is a subbundle of the tangent bundle of $P$, and $P$ is a submanifold
of $M$, we consider $\widetilde{v}$ as a section of $TP$ and $\widetilde{
w}$ as a section of $TM_{|P}$. Restricting the hyper-K\"{a}hler connection $
\nabla $ onto $P$, we consider 
$\nabla _{\widetilde{v}}\widetilde{w}$ as a section
of $TM_{|P}$. Recall that there is a direct sum decomposition 
\begin{equation}
TM_{|P}=U\oplus \mathfrak{g\oplus }I_{1}\mathfrak{g\oplus }I_{2}\mathfrak{g\oplus }I_{3}
\mathfrak{g}.
\end{equation}
Let $\theta $ be the projection from $TM_{|P}$ onto its direct summand $U$.
Since $\mathfrak{g}$ is orthogonal to the distribution $U$, and $U$ is
hypercomplex invariant, $\theta $ is an orthogonal projection. Define 
\begin{equation}
D_{v}w:=d\pi (\theta (\nabla _{\widetilde{v}}\widetilde{w})).\quad 
\mbox{i.e. \quad }\widetilde{D_{v}w}=\theta (\nabla _{\widetilde{v}}
\widetilde{w}).
\end{equation}
Now we have to prove that it is a HKT-connection. 

We claim that the connection $D$ preserves the hypercomplex
structure. This claim is equivalent to $D_{v}(I_aw)=I_aD_{v}w$. Lifting
to $U$, it is equivalent to 
$\theta (\nabla _{\widetilde{v}}I_a\widetilde{w})=I_a\theta (\nabla _{
\widetilde{v}}\widetilde{w})$.
Since the direct sum decomposition is invariant of the hypercomplex
structure, the projection map $\theta $ is hypercomplex. Therefore, it
commutes with the complex structures. Then the above identity is equivalent
to 
$\theta (\nabla _{\widetilde{v}}I_a\widetilde{w})=\theta (I_a\nabla _{
\widetilde{v}}\widetilde{w})$.
This identity holds because $\nabla $ is hypercomplex. 

To verify that the connection $D$ preserves the Riemannian metric $h$, 
let $u,v,$ and $w$ be vector fields on $N$. The identity
$uh(v,w)-h(D_{u}v,w)-h(v,D_{u}w)=0$
is equivalent to the following identity on $P$:
$
\widetilde{u}g(\widetilde{v},\widetilde{w})-g(\theta (\nabla _{\widetilde{u}}
\widetilde{v}),\widetilde{w})-g(\widetilde{v},\theta (\nabla _{\widetilde{u}}
\widetilde{w}))=0.
$
Since $\theta $ is the orthogonal projection along $\mathfrak{g}$, the above
identity is equivalent to 
$
\widetilde{u}g(\widetilde{v},\widetilde{w})-g(\nabla _{\widetilde{u}}
\widetilde{v},\widetilde{w})-g(\widetilde{v},\nabla _{\widetilde{u}}
\widetilde{w})=0.
$
This identity on $P$ is satisfied because $\nabla $ is a HKT-
connection. 

Finally, we have to verify that the torsion of the
connection $D$ is totally skew-symmetric. By definition and
the fact that $\theta$ is an orthogonal projection, the torsion of 
$D$ is
$T^{D}(u,v,w)=
g(\nabla _{\widetilde{u}}\widetilde{v},\widetilde{w})-g(\nabla _{
\widetilde{v}}\widetilde{u},\widetilde{w})-g(\widetilde{[u,v]},\widetilde{w})$.
Note that $[\widetilde{u},\widetilde{v}]$ is a vector tangent to $P$ such
that $d\pi \circ \theta ([\widetilde{u},\widetilde{v}])=[d\pi (\widetilde{u}
),d\pi (\widetilde{v})]=[u,v].$ Therefore, $[\widetilde{u},\widetilde{v}]$
and $\widetilde{[u,v]}$ differ by a vector in $\mathfrak{g}$. Since the Killing
vector fields are orthogonal to the hypercomplex distribution, 
$g(\widetilde{[u,v]},\widetilde{w})=g([\widetilde{u},\widetilde{v}],
\widetilde{w})$.
Then we have 
$T^{D}(u,v,w)=T^{\nabla }(\widetilde{u},\widetilde{v},\widetilde{w})$.
This is totally skew-symmetric because the connection $\nabla $ is the
Bismut connection on $M$. \eproof

Suppose that the group $G$ is one-dimensional. Let $X$ be
the Killing vector field generated by $G$. The hypercomplex distribution $U$
and the horizontal distribution $H$ are identical if and only if the
1-forms $I_{1}d\nu _{1}=I_{2}d\nu _{2}=I_{3}d\nu _{3}$ are pointwisely
proportional to the 1-form $\iota _{X}g$ along the level set $P$. i.e. for
any tangent vector $Y $ to $P$, 
$I_ad\nu_a(Y) =fg(X,Y)$. Equivalent 
$d\nu_a =f\iota _{X}F_a$.
In the next example, we shall make use of this observation.

\subsection{Example: HKT-Structure on 
$\mathcal{V}\left( \cp^2\right) =S^{1}\times (SU(3)/U(1))$}

We construct a HKT-structure on $\mathcal{V}\left( \cp^2\right) 
$by a $U(1)$-reduction from a HKT-structure on 
$\bH^{3}\backslash\{0\}.$ 
Choose a hypercomplex structure on 
$\bR^{6}\cong \bC^{3}\oplus \bC^{3}$ by 
\begin{equation}
I_{1}(\chi ,\varrho )=(i\chi ,-i\varrho ),\quad I_{2}(\chi ,\varrho
)=(i\varrho ,i\chi ),\quad I_{3}(\chi ,\varrho )=(-\varrho ,\chi ).
\end{equation}
It is apparent that the holomorphic coordinates with these complex
structures are $(\chi ,\overline{\varrho })$, $(\chi +\varrho ,\overline{
\chi }-\overline{\varrho }),$ and $(\varrho -i\chi ,\overline{\varrho }-i
\overline{\chi })$ respectively. 

As in Example \ref{hopf}, the hyper-K\"ahler potential for
 the Euclidean metric $g$ on 
$(\bC^{3}\oplus \bC^{3})\backslash\{0\}$
is $\mu=\frac12(|\chi|^2+|\varrho|^2)$. We apply Proposition
\ref{modification} to $f(\mu)=\ln\mu$ to obtain a new
HKT-metric 
\begin{equation}
{\hat g}=\frac{1}{\mu}g-\frac{1}{\mu^2}(d\mu\otimes
d\mu+I_1d\mu\otimes I_1d\mu+ I_2d\mu\otimes I_2d\mu + I_3d\mu\otimes
I_3d\mu).
\end{equation}

Define a hypercomplex moment map 
$\nu=(\nu_{1},\nu _{2},\nu _{3})$ by 
\begin{equation}
\nu _{1}(\chi ,\varrho )=|\chi |^{2}-|\varrho |^{2},\quad (\nu _{2}+i\nu
_{3})(\chi ,\varrho )=2\left\langle \chi ,\varrho \right\rangle .
\end{equation}
where $\left\langle ,\right\rangle $ is a Hermitian inner product on $\bC
^{3}$. Let $\Gamma \cong \U(1)$ be the one-parameter group acting on $(\bC
^{3}\oplus \bC^{3})\backslash\{0\}$ defined by 
\begin{equation}
(t;(\chi ,\varrho ))\mapsto (e^{it}\chi ,e^{it}\varrho ).
\end{equation}
Let $\left\langle r\right\rangle $ be the integer group generated by a real
number between $0$ and $1$. It acts on 
$(\bC^{3}\oplus \bC^{3})\backslash\{0\}$ by 
\begin{equation}
(n;(\chi ,\varrho ))\mapsto (r^{n}\chi ,r^{n}\varrho ).
\end{equation}
Both $\Gamma $ and $\left\langle r\right\rangle $ are groups of hypercomplex
automorphisms leaving the zero level set of $\nu $ invariant. Then the
quotient space $\nu ^{-1}(0)/\Gamma $ is a hypercomplex reduction. The
discrete quotient space 
$\mathcal{V}=\nu ^{-1}(0)/\Gamma\times\left\langle r\right\rangle$ 
is a compact hypercomplex manifold. From the homogeneity of the metric $\hat g$,
we see that both $\Gamma$ and the discrete group $\left\langle r\right\rangle$
are group of isometries for the metric $\hat g$. Therefore, the
quotient space $\mathcal{V}$ inherits a hyper-Hermitian metric.

On $(\bC^{3}\oplus \bC^{3})\backslash\{0\}$, 
the real vector field generated by the group $\Gamma $ is 
\[
X=i\chi \frac{\partial }{\partial \chi }-i\overline{\chi }\frac{\partial }{
\partial \overline{\chi }}-i\overline{\varrho }\frac{\partial }{\partial 
\overline{\varrho }}+i\varrho \frac{\partial }{\partial \varrho }.
\]
Let ${\hat F}_a$ be the K\"ahler form for the HKT-metric $\hat g$. We check
that
$d\nu_a=-2\mu \iota_X{\hat F}_a$. Therefore,
Theorem \ref{reduction} implies that the quotient space $\cal V$
inherits a HKT-structure.

 Note that if 
$(\chi ,\varrho )$ is a point in the zero level set, then it represents a
pair of orthogonal vectors. Therefore, the triple 
$(\frac{\chi }{|\chi |},
\frac{\varrho }{|\varrho |},
\frac{\overline{\chi }}{|\chi |}\times 
\frac{\overline{\varrho }}{|\varrho |})$ forms an element in the
matrix group $\SU(3).$ The action of $\Gamma $ induces an action on $\U(3)$ by
the left multiplication of $\diag(e^{it},e^{it},e^{-2it})$. 
Denote the 
$\Gamma $-coset of 
$(\frac{\chi }{|\chi |},
\frac{\varrho }{|\varrho |},
\frac{\overline{\chi }}{|\chi |}\times 
\frac{\overline{\varrho }}{|\varrho |})$
 by $[\frac{\chi }{|\chi |},
\frac{\varrho }{|\varrho |},
\frac{\overline{\chi }}{|\chi |}\times 
\frac{\overline{\varrho }}{|\varrho |}]$. The quotient space $\mathcal{V}$ is
isomorphic to the product space $S^{1}\times \SU(3)/\U(1)$. The quotient map is 
\[
(\chi ,\varrho )\mapsto \left( \exp \left( 2\pi i\frac{\ln |\chi |}{\ln r}
\right) ,\quad \lbrack \frac{\chi }{|\chi |},\frac{\varrho }{|\varrho |},
\frac{\overline{\chi }}{|\chi |}\times 
\frac{\overline{\varrho }}{|\varrho |}]
\right) .
\]

\noindent{\bf Remark:} A fundamental question on HKT-structures remains open.
Does every hypercomplex
manifold admit a metric such that it is a HKT-structure?

\vspace{.1in} 

\noindent{\bf   Acknowledgment} We thank G.W. Gibbons for
introducing the topic in this paper to us. The second author thanks J.-P.
Bourguignon for providing an excellent research environment at the I.H.E.S..

\end{document}